\documentclass{article}

\usepackage[T1]{fontenc}
\usepackage[utf8]{inputenc}
\usepackage[english]{babel}
\usepackage{newtxtext}

\usepackage{amsmath}
\usepackage{amssymb}
\usepackage{bm}
\usepackage{bbm}
\usepackage{spalign}
\usepackage{enumerate}
\usepackage{fix-cm}

\usepackage{caption}
\usepackage{subcaption}

\usepackage[obeyspaces]{url}

\usepackage{graphicx}
\usepackage{float}

\usepackage{fancyvrb}

\usepackage[theme=grayscale]{jlcode}

\usepackage[section]{algorithm}
\usepackage{algorithmic}

\usepackage{xcolor}

\usepackage[colorinlistoftodos]{todonotes}

\usepackage{booktabs}

\title{\href{https://github.com/ZIB-IOL/FrankWolfe.jl}{\texttt{FrankWolfe.jl}}: A high-performance and flexible toolbox \\ for Frank-Wolfe algorithms and Conditional Gradients}

\author{
       \name Mathieu Besançon \email \href{mailto:besancon@zib.de}{besancon@zib.de} \\
       \addr Zuse Institute Berlin, Germany
       \AND
\name Alejandro Carderera \email  \href{mailto:alejandro.carderera@gatech.edu}{alejandro.carderera@gatech.edu}\\
       \addr Georgia Institute of Technology, USA\\
       Zuse Institute Berlin, Germany 
       \AND
       \name Sebastian Pokutta \email \href{mailto:pokutta@math.tu-berlin.de}{pokutta@math.tu-berlin.de} \\
       \addr Technische Universit\"at Berlin, Germany\\
       Zuse Institute Berlin, Germany
}

\DeclareMathOperator*{\argmin}{argmin}
\newcommand{\innp}[2]{\left\langle #1, #2 \right\rangle}
\newcommand{\norm}[1]{\left\| #1 \right\|}
\newcommand{\vx}{\mathbf{x}}
\newcommand{\vy}{\mathbf{y}}
\newcommand{\vd}{\mathbf{d}}
\newcommand{\vc}{\mathbf{c}}

\usepackage{jmlr2e}
\usepackage[capitalize,nameinlink]{cleveref}[0.21]
\usepackage{natbib}
\setcitestyle{authoryear,round,citesep={;},aysep={,},yysep={;}}
\usepackage[vvarbb]{newtxmath}
\usepackage{multirow}

\begin{document}

\maketitle

\begin{abstract}
We present \texttt{FrankWolfe.jl}, an open-source implementation of several popular Frank-Wolfe and Conditional Gradients variants for first-order constrained optimization.
The package is designed with flexibility and high-performance in mind, allowing for easy extension and relying on few assumptions regarding the user-provided functions.
It supports Julia's unique multiple dispatch feature, and interfaces smoothly with generic linear optimization formulations using \texttt{MathOptInterface.jl}.
\end{abstract}

\section{Introduction}

We provide an open-source software package released under the MIT license, for \emph{Frank-Wolfe (FW)} and
\emph{Conditional Gradients} algorithms implemented in Julia \citep{bezanson2017julia}.
Its focus is a general class of constrained convex optimization problems of the form:
\begin{align*}
\min_{\vx \in \mathcal{C}} f(\vx)
\end{align*}
where $\mathcal{X}$ is a Hilbert space,
$\mathcal{C} \subseteq \mathcal{X}$ is a compact convex set, and
$f: \mathcal{C} \to \mathbb{R}$ is a convex, continuously differentiable function.

Although the Frank-Wolfe algorithm and its variants have been studied for more than half a century and have gained a lot of
attention for their theoretical and computational properties, no general-purpose off-the-shelf implementation exists.
The purpose of the package is to become a reference open-source implementation
for practitioners in need of a flexible and efficient first-order method and for researchers developing and
comparing new approaches on similar classes of problems.

\section{Frank-Wolfe algorithms}

The Frank-Wolfe algorithm \citep{frank1956algorithm} (also known as the Conditional Gradient algorithm \citep{polyak66cg}),
is a first-order algorithm for constrained optimization that avoids the use of projections at each iteration.
For the sake of exposition, we confine ourselves to the Euclidean setup. The main ingredients that the FW algorithm leverages are:
\begin{enumerate}
\item \textbf{First-Order Oracle} (FOO): Given $\vx \in \mathcal{C}$, the oracle returns $\nabla f(\vx)$.
\item \textbf{Linear Minimization Oracle} (LMO): Given $\vd \in \mathcal{X}$, the oracle returns:
\begin{align} \label{prob:lmo}
    \mathbf{v} \in \argmin_{\vx \in \mathcal{C}} \langle \vd, \vx\rangle
\end{align}
\end{enumerate}
The simplest version of the algorithm (shown in \cref{vanillafw}) builds a linear approximation to the function at a
given iterate, using first-order information, and minimizes this approximation over the feasible region (Line~\ref{fw_vertex} in \cref{vanillafw}).
A schematic representation of the step is described in \cref{fig:schematic_CG_step} where the {\color{blue}blue} curves represent the contour lines of $f(\vx)$, and the {\color{red}red} lines represent the contour lines of the linear approximation built at $\vx_t$.
The new iterate is then computed as a convex combination of the current iterate and the linear minimizer $\mathbf{v}_t$ from the LMO (Line~\ref{convex_combination} in \cref{vanillafw}).
Since $\mathbf{v}_t$ is an extreme point of the feasible region, the new iterate is feasible by the convexity of $\mathcal{C}$ and remains so throughout the algorithm.
Alternatively, one can view each iteration of the Frank-Wolfe algorithm as finding the direction that is best aligned with the negative of the gradient using only the current iterate and the extreme points of the feasible region.

\begin{minipage}{0.40\textwidth}
\begin{algorithm}[H]
  \caption{Frank-Wolfe algorithm}
\label{vanillafw}
\begin{algorithmic}[1]
  \REQUIRE Point $\vx_0\in \mathcal{C}$, function \(f\),
  stepsizes $\gamma_t >0$.
\ENSURE Iterates $\vx_1, \dotsc \in \mathcal{C}$.
\FOR{$t=0$ \TO \dots}
  \STATE$\mathbf{v}_t\leftarrow\argmin_{\mathbf{v}\in\mathcal{C}}\innp{\nabla f(\vx_t)}{\mathbf{v}}$
    \label{fw_vertex}
\STATE$\vx_{t+1}\leftarrow \vx_t+\gamma_t(\mathbf{v}_t-\vx_t)$\label{convex_combination}
\ENDFOR
\end{algorithmic}
\end{algorithm}
\end{minipage}
\begin{minipage}{0.60\textwidth}
\begin{figure}[H]
\centering
  \includegraphics[width=0.7\linewidth]{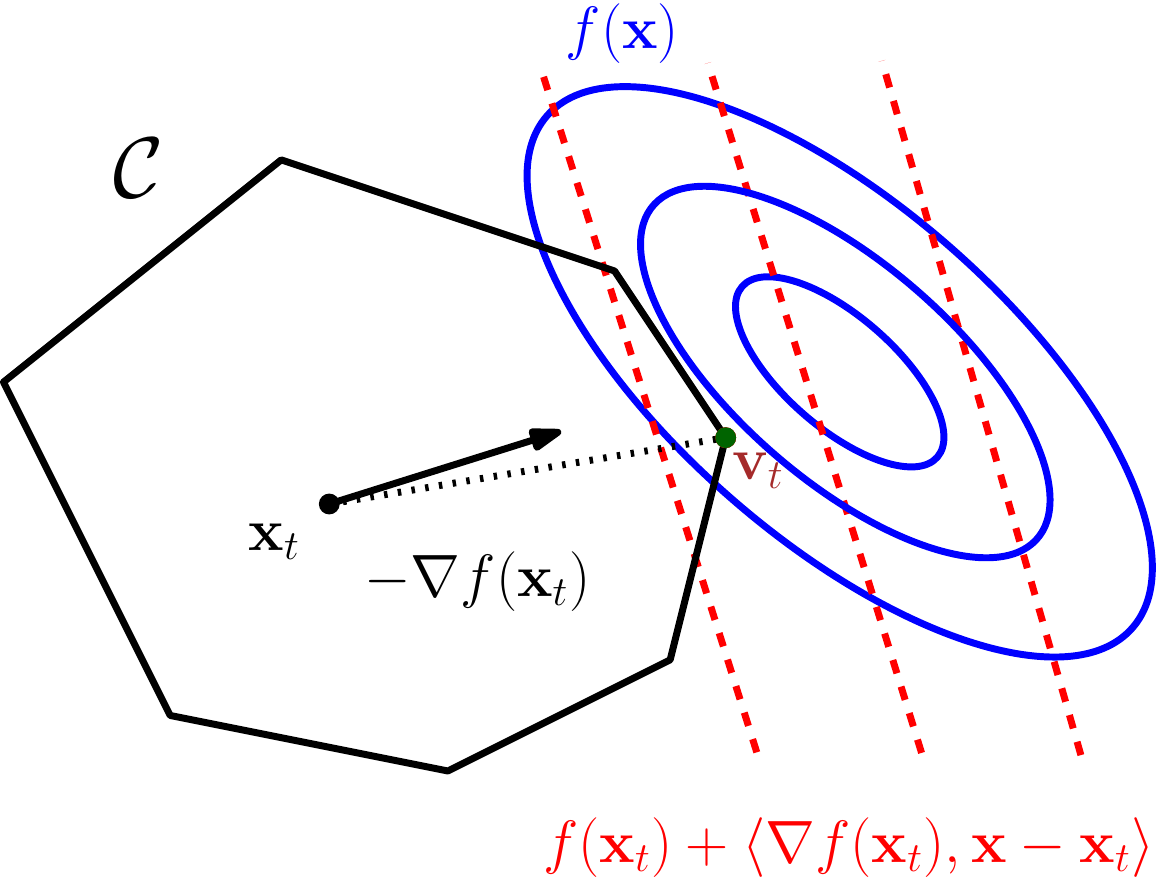}
 \caption{Schematic of a Frank-Wolfe step at iteration $t$.}
    \label{fig:schematic_CG_step}
\end{figure}
\end{minipage} \medskip

The Frank-Wolfe algorithm and its variants have recently gained a lot of attention in the Machine Learning community due to their attractive properties. For example, as the iterates are formed by convex combinations of extreme points of the feasible set, the solutions outputted by the algorithm can often be expressed as sparse combinations of a few extreme points. As first-order methods, they only assume
mild requirements with respect to accessing the objective function, and only require being able to compute the gradient of $f$, and not its Hessian.
Furthermore, these algorithms only require being able to compute the minimum of a linear function over the feasible region.
This is particularly advantageous in applications where computing a projection (in essence, solving a quadratic problem)
is much more computationally expensive than solving a linear optimization problem \citep{combettes21complexity}.
For example, if $\mathcal{C}$ is the set of matrices in $\mathbb{R}^{n \times m}$ of bounded nuclear norm, computing the projection of a matrix onto $\mathcal{C}$ requires computing a full Singular Value Decomposition.
On the other hand, solving a linear optimization problem over $\mathcal{C}$ only requires computing the left and right singular vectors associated with the largest singular value.

Another interesting property of the FW algorithm is that when minimizing a convex function $f$,
the primal gap $f(\vx) - f(\vx^*)$ is bounded by a quantity dubbed the \emph{Frank-Wolfe gap} or the \emph{dual gap} defined as $\innp{\nabla f(\vx_t)}{\vx_t - \mathbf{v}_t}$.
This follows directly from the convexity of $f$ and the maximality of $\mathbf{v}_t$ since:
\begin{align*}
	 f(\vx) - f(\vx^*) \leq \innp{\nabla f(\vx_t)}{\vx_t - \vx^*} \leq \max_{\mathbf{v} \in \mathcal{C}}\innp{\nabla f(\vx_t)}{\vx_t - \mathbf{v}} = \innp{\nabla f(\vx_t)}{\vx_t - \mathbf{v}_t}.
\end{align*}
The dual gap (and its stronger variants) is extremely useful as a measure of progress and is often used as a stopping criterion when running Frank-Wolfe algorithms.

\subsection{Algorithm variants}

The basic algorithm presented in \cref{vanillafw} requires few assumptions on
the problem structure, making it a versatile baseline option.
It is however not able to attain linear convergence in general, even when the objective function is strongly convex,
which has led to the development of algorithmic variants that enhance the performance of the original algorithm,
while maintaining many of its advantageous properties.
We summarize below the central ideas of the variants implemented in the package and highlight in \cref{table:zero} key properties that can drive the choice of a variant on a given use case.
We remark that most variants also work for the case where the objective is nonconvex, providing some locally optimal solution in this case.
The \textbf{Sparsity} column refers to the number of extreme points of $\mathcal{C}$ composing the iterates, as customary for FW methods, and not to the number of non-zero entries
of the iterate. Some algorithms are more robust to ill-conditioned problems, for instance when some operations are rendered challenging by numerical precision, this aspect is captured in the \textbf{Numerical Stability} column.

\paragraph{Standard Frank-Wolfe} The simplest Frank-Wolfe variant is presented in \cref{vanillafw}. It has the lowest memory requirements out of all the variants, as in its simplest form only requires keeping track of the current iterate. As such it is suited for extremely large problems.
However, in certain cases, this comes at the cost of speed of convergence in terms of iteration count, when compared to other variants. As an example, when minimizing a strongly convex and smooth function over a polytope this algorithm might converge sublinearly, whereas the three variants that will be presented next converge linearly.

\begin{figure*}[th!]
\centering
\hspace{\fill}
\subfloat[Standard Frank-Wolfe]{{\includegraphics[width=6cm]{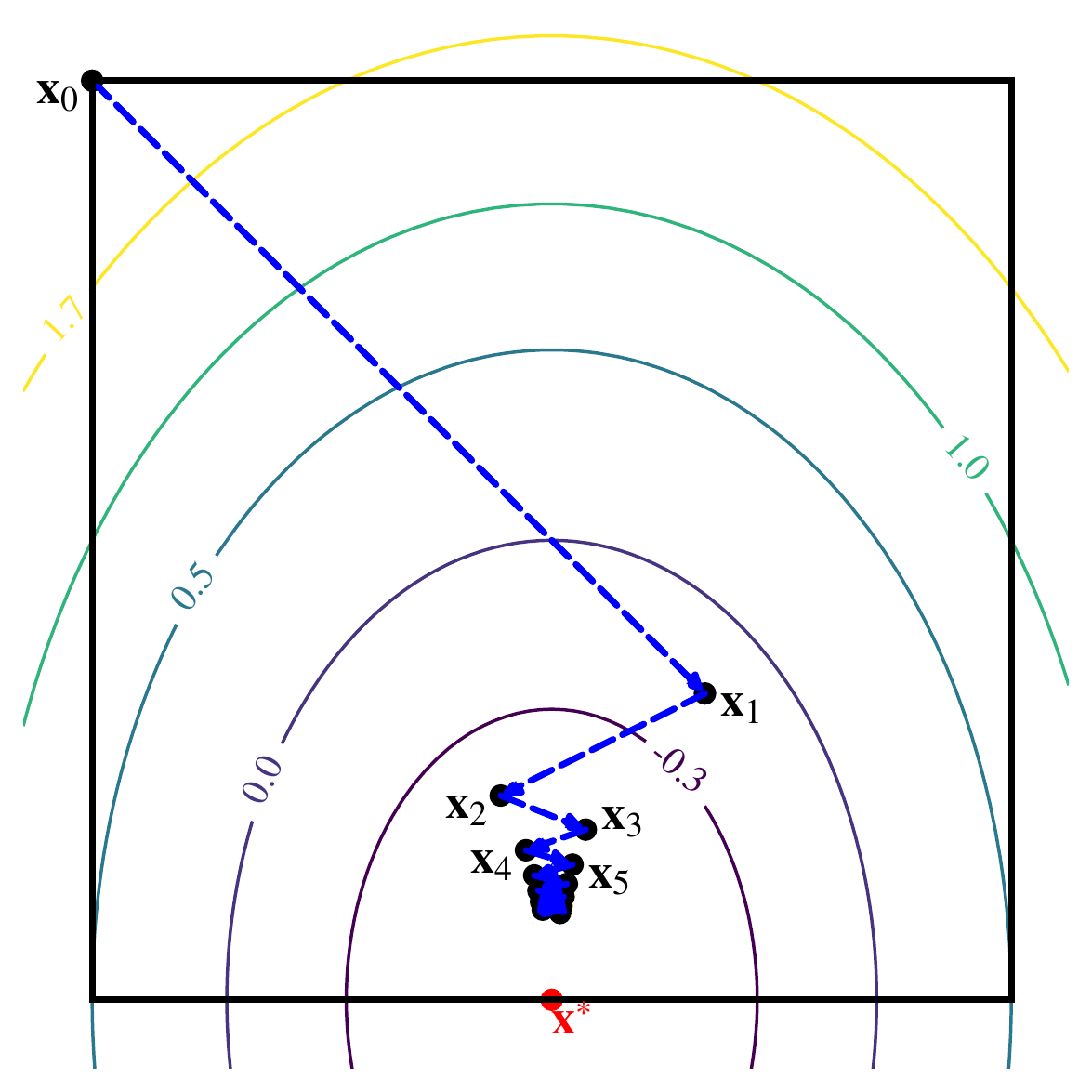} }\label{fig:Vanilla_convergence}}%
\hspace{\fill}
\subfloat[Away-step Frank-Wolfe]{{\includegraphics[width=6cm]{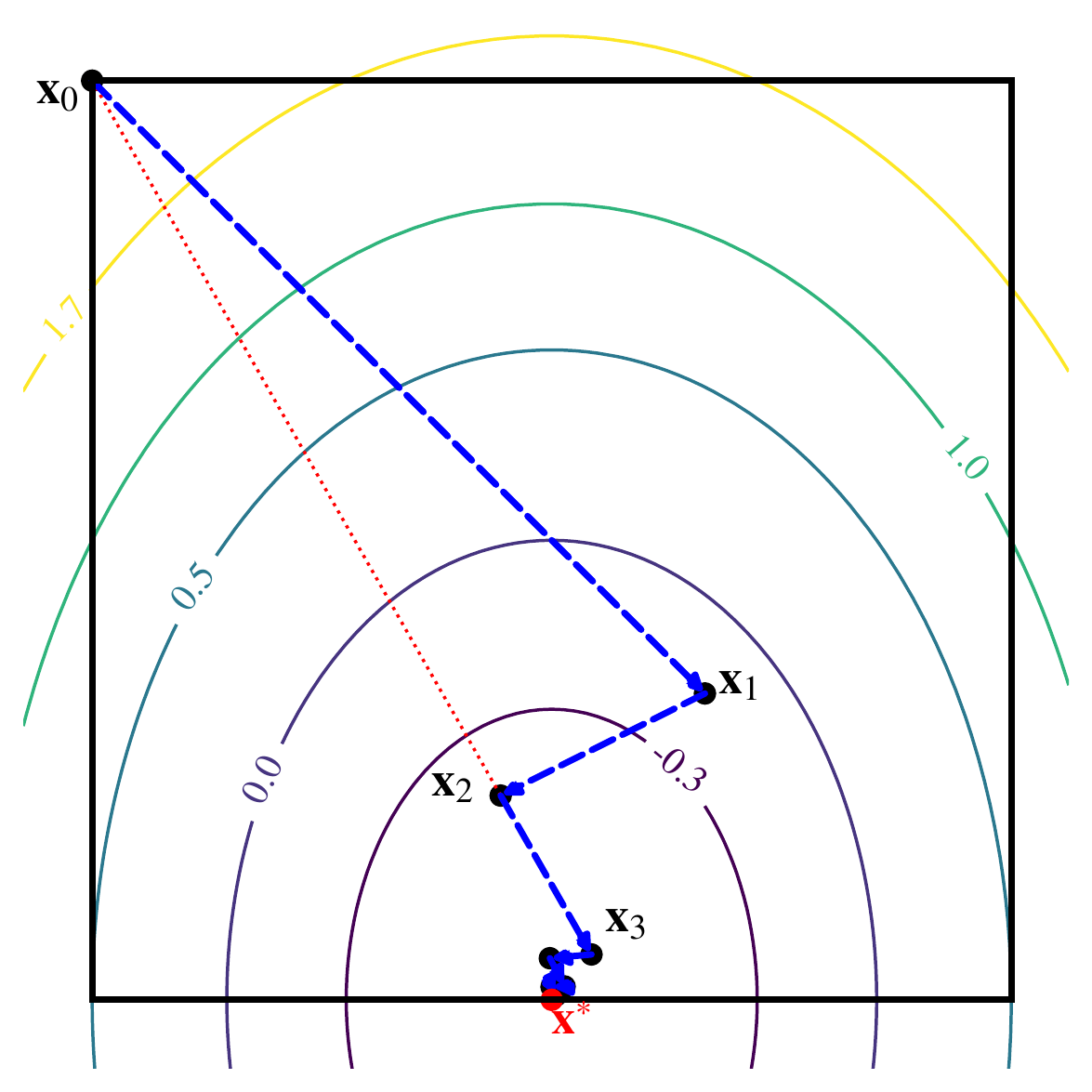} }\label{fig:Away_convergence}}%
\hspace{\fill}
\caption{
Convergence comparison of the standard Frank-Wolfe algorithm and the Away-step Frank-Wolfe algorithm when minimizing a quadratic function (contour lines depicted) over a simple polytope using exact line search.
Note that the minimizer is on a face of the polytope, which causes the Frank-Wolfe algorithm to slowly \emph{zig-zag} towards the solution, causing sublinear convergence.
At the iterate $\vx_2$, AFW performs an away-step, moving away from the extreme point $\vx_0$ that belongs to the active set, the set of vertices that compose the convex decomposition of the current iterate. The weight of $\vx_0$ in the current decomposition of the iterate is thus reduced. The alignment of the away direction with the gradient is better than the alignment that would be obtained by a standard Frank-Wolfe step.
}
\label{fig:convergence_comparison}
\end{figure*}

\paragraph{Away-step Frank-Wolfe} The most popular among the Frank-Wolfe variants is the \emph{Away-step Frank-Wolfe} (AFW) algorithm \citep{guelat1986some, lacoste15}. While the FW algorithm presented in \cref{vanillafw} can only move \emph{towards} extreme points of $\mathcal{C}$, the AFW can move \emph{away} from some extreme points of $\mathcal{C}$, hence the name of the algorithm. To be more specific, the AFW algorithm moves away from vertices in its \emph{active set} at iteration  $t$, denoted by $\mathcal{S}_t$, which contains the set of vertices $\mathbf{v}_k$ for $k<t$ that allow us to recover the current iterate as a convex combination. This algorithm expands the range of directions that the FW algorithm can move along, at the expense of having to explicitly maintain the current iterate as a convex decomposition of extreme points. See \cref{fig:convergence_comparison} for a schematic of the behavior of the two algorithms for a simple example.

\paragraph{Lazifying Frank-Wolfe variants} One running assumption for the two previous variants is that calling the LMO is cheap. There are many applications where calling the LMO in absolute terms is costly (but is cheap in relative terms when compared to performing a projection). In such cases, one can attempt to \emph{lazify} FW algorithms, 
to avoid having to compute $\argmin_{\mathbf{v}\in\mathcal{C}}\innp{\nabla f(\vx_t)}{\mathbf{v}}$ by calling the LMO, settling for finding approximate minimizers that guarantee enough progress \citep{pok17lazy}. This allows us to substitute the LMO by a \emph{Weak Separation Oracle} while maintaining essentially the same convergence rates. In practice, these algorithms search for an appropriate extreme point among the extreme points in a cache, or the extreme points in the active set $\mathcal{S}_t$, and can be much faster in wall-clock time. In the package, both AFW and FW have lazy variants while the \emph{Blended Conditional Gradient} (BCG) algorithm is lazified by design.

\paragraph{Blended Conditional Gradients} The FW and AFW algorithms, and their lazy variants share one feature: they attempt to make primal progress over a reduced set of extreme points.
The AFW algorithm does this through away steps (which do not increase the cardinality of the active set), and the lazy variants do this through the use of previously exploited extreme points.
A third strategy that one can follow is to explicitly \emph{blend} Frank-Wolfe steps with gradient descent steps over the convex hull of the active set (note that this can be done without requiring a projection oracle over $\mathcal{C}$, thus making the algorithm projection-free). This results in the \emph{Blended Conditional Gradient} (BCG) algorithm \citep{pok18bcg}, which attempts to make as much progress as possible through the convex hull of the current active set $\mathcal{S}_t$ until it automatically detects that in order to further make further progress it requires additional calls to the LMO.

\paragraph{Stochastic Frank-Wolfe}
In many problem instances, evaluating the FOO at a given point is prohibitively expensive. In some of these cases, one can leverage a \emph{Stochastic First-Order Oracle} (SFOO), from which one can build a gradient estimator \citep{stoch51}. This idea, which has powered much of the success of deep learning, can also be applied to the Frank-Wolfe algorithm \citep{hazan2016variance}, resulting in the \emph{Stochastic Frank-Wolfe} (SFW) algorithm and its variants.
We provide a flexible interface with the possibility of providing rules for both the batch size and the momentum as series indexed by the iteration count.
This interface allows users to implement the growing batch size necessary for the stochastic variant of \cite{hazan2016variance} and the fixed-batch damped momentum of \cite{mokhtari2020stochastic}.

\begin{table}[H]
\centering
\begin{tabular}{ccccccc}
\toprule
\multirow{2}{*}{\textbf{Algorithm}} & \multicolumn{2}{c}{\textbf{Convergence}} & \multirow{2}{*}{\textbf{Sparsity}} & \textbf{Numerical} & \multirow{2}{*}{\textbf{Active Set?}} \\
& Progress/iteration   & Time/iteration   & & \textbf{Stability} & \\
\midrule
FW & Low & Low & Low & High & No \\
AFW & Medium & Medium-High & Medium & Medium-High & Yes  \\
SFW & Low & Low & Low & High & No  \\
\hline
Lazy FW & Low & Low & High & High & Yes  \\
Lazy AFW & Medium & Medium & High & Medium-High & Yes  \\
BCG & Medium-High & Medium & High & Medium & Yes  \\
\bottomrule
\end{tabular}
 \caption{Comparison of the characteristics of the algorithms in the package, when applied to a general problem. The algorithms below the horizontal line in the middle of the table are \emph{lazy} by design, while the algorithms above the line are not. In practice, lazy algorithms lead to a higher level of sparsity than non-lazy algorithms, moreover they tend also tend to run faster than their non-lazy counterparts when the LMO is computationally expensive, as they try to avoid calls to the LMO when an appropriate vertex already exists in the active set.}
 \label{table:zero}
\end{table}

\section{The \texttt{FrankWolfe.jl} package}

The package offers the first implementation of Frank-Wolfe variants in Julia and more broadly
tackles constrained optimization problems in a way that complements well the thriving
Julia optimization ecosystem.\\

Unlike disciplined convex frameworks or algebraic modeling languages such as
\texttt{Convex.jl} \citep{udell2014convex} or \texttt{JuMP.jl} \citep{dunning2017jump,legat2020mathoptinterface},
our framework allows for arbitrary Julia functions defined outside of a Domain-Specific Language.
Users can provide their gradient implementation or leverage one of the many automatic differentiation
packages available in Julia.

Several ecosystems have emerged for non-linear optimization, but these often only handle specific types of constraints.
Most notably, \texttt{JSOSolvers.jl} \citep{orban-siqueira-jsosolvers-2020} offers a
trust-region second-order method for bound constraints. On the other hand,
\texttt{Optim.jl} \citep{Optim.jl-2018} implements an interior-point method handling bound and non-linear constraints.
Constrained first-order methods based on proximal operators are implemented in \texttt{StructuredOptimization.jl}
and \texttt{ProximalAlgorithms.jl} \citep{antonello2018proximal} but only allow specific functions defined through the exposed
Domain-Specific Language.

One central design principle of \texttt{FrankWolfe.jl} is to rely on few assumptions
regarding the user-provided functions, the atoms returned by the LMO, and
their implementation. The package works for instance out of the box when the LMO
returns Julia subtypes of \jlinl{AbstractArray}, representing finite-dimensional
vectors, matrices or higher-order arrays.

Another design principle has been to favor in-place operations and reduce memory allocations when possible, since these
can become expensive when repeated at all iterations.
This is reflected in the memory emphasis mode (the default mode for all algorithms),
in which as many computations as possible are performed in-place, as well as in the gradient
interface, where the gradient function is provided with a storage variable to write into
in order to avoid reallocating a new variable every time a gradient is computed: \jlinl{grad!(storage, x)}.

The performance difference due to this additional array allocation can be quite pronounced for problems in large dimensions.
For example, computing an out-of-place gradient of size $7.5$GB on a state-of-the-art machine is about $8$ times slower than an in-place update
due to the additional allocation of the large gradient object.

Finally, default parameters are chosen to make all algorithms as robust as possible out of the box, while allowing extension and fine tuning for advanced users.
For example, the default step size strategy for all (but the stochastic variant) is the adaptive step size rule of \citet{pedregosa2018step}, which in computations not only usually outperforms both line search and the short step rule by dynamically estimating the Lipschitz constant but also overcomes several issues with the limited additive accuracy of traditional line search rules.
Similarly, the BCG variant automatically upgrades the numerical precision for certain subroutines if numerical instabilities are detected.

\subsection{Linear minimization oracle interface}

One key step of FW algorithms is the linear minimization step which, given first-order information
at the current iterate, returns an extreme point of the feasible region that minimizes the linear approximation
of the function. It is defined in \texttt{FrankWolfe.jl} using a single function:

\begin{jllisting}
function compute_extreme_point(lmo::LMO, direction::D; kwargs...)::V
  # ...
end
\end{jllisting}

The first argument \jlinl{lmo} represents the linear minimization oracle for the specific problem.
It encodes the feasible region $\mathcal{C}$, but also some algorithmic parameters or state.
This is especially useful for the lazified FW variants, as in these cases the LMO types can
take advantage of caching, by storing the extreme 
points that have been computed in previous iterations and then looking up extreme points 
from the cache before computing a new one. 

The package implements LMOs for commonly encountered feasible regions including $L_p$-norm balls,
$K$-sparse polytopes, the Birkhoff polytope, and the nuclear norm ball for matrix spaces, leveraging known
closed forms of extreme points. The multiple dispatch mechanism allows for different implementations
of a single LMO with multiple direction types. The type \jlinl{V} used to represent the computed
extreme point is also specialized to leverage the properties of extreme points of the feasible region.
For instance, although the Birkhoff polytope is the convex hull of all doubly stochastic matrices of a given
dimension, its extreme points are permutation matrices that are much sparser in nature.
We also leverage sparsity outside of the traditional sense of nonzero entries.
When the feasible region is the nuclear norm ball in $\mathbb{R}^{N\times M}$, the extreme points are rank-one matrices.
Even though these extreme points are dense, they can be represented as the outer product of two vectors
and thus be stored with $\mathcal{O}(N+M)$ entries instead of $\mathcal{O}(N\times M)$ for the equivalent
dense matrix representation. The Julia abstract matrix representation allows the user and the library to interact
with these rank-one matrices with the same API as standard dense and sparse matrices.

In some cases, users may want to define a custom feasible region that does not admit a closed-form
linear minimization solution. We implement a generic LMO based on \texttt{MathOptInterface.jl} (MOI) \citep{legat2020mathoptinterface},
thus allowing users on the one hand to select any off-the-shelf LP, MILP, or conic solver suitable for their problem,
and on the other hand to formulate the constraints of the feasible domain using the \texttt{JuMP.jl} or \texttt{Convex.jl} DSL.
Furthermore, the interface is naturally extensible by users who can define their own LMO and implement the corresponding \jlinl{compute_extreme_point} method.

\subsection{Numeric type genericity}

The package was designed from the start to be generic over both the used numeric types and data structures.
Numeric type genericity allows running the algorithms in extended fixed or arbitrary precision, e.g., the package works out-of-the-box with \jlinl{Double64} and \jlinl{BigFloat} types.
Extended precision is essential for high-dimensional problems where the condition number of the gradient computation
becomes too high.
For some well-conditioned problems, reduced precision is sometimes sufficient to achieve the desired tolerance.
Furthermore, it opens the possibility of gradient computation and LMO steps on hardware accelerators such as GPUs.

\section{Application examples}

We will now present a few application examples that highlight specific features of the package.
The full code of each example (and several more) can be found in the examples folder of the repository.
The progress of different algorithms against time and number of iterations is displayed to illustrate their
properties and strengths summarized in \cref{table:zero} in various settings.

\subsection{Polynomial Regression}

Data encountered in many real-world applications can
be modeled as coming from data-generating processes, which non-linearly 
map a set of input features to an output space. This non-linear mapping can be often be 
modeled with a relatively sparse 
linear combination of simpler non-linear basis functions. 
For example, if we denote the scalar output by $y \in \mathbb{R}$, the input features
by $\vx \in \mathbb{R}^n$, and the library of basis functions by $f_i \colon \mathbb{R}^n \to \mathbb{R}$ with $i\in [1, m]$ we might have that the data is generated as:
\begin{align*}
	y = \sum_{i = 1}^m c_i f_i(\vx),
\end{align*}
where $c_i \in \mathbb{R}$, and only a few 
of the basis functions participate in the data generating process, i.e., many of
the coefficients satisfy $c_i = 0$. Our task is to recover the non-zero $c_i$ coefficients 
when given input/output measurements $\{\vx_j, y_j\}_{j = 1}^N$. In the absence of noise, and collecting the data points 
into vectors $\vy\in \mathbb{R}^N$, $\vc\in \mathbb{R}^m$, and an appropriate matrix $A \in \mathbb{R}^{N\times m}$ 
we can write $\vy = A \vc$. 
We focus on the case 
where $N \ll m$, where we have more unknowns then equations, inside the 
realm of what is known as compressed sensing (see \cite{candes2008introduction} for an 
overview), where we need to impose 
additional structure on our learning problem in 
order for it to make sense. 
This is where the sparsity assumption comes into play, 
as we know that only a few of the 
$c_i$ coefficients are non-zero. Due to the fact that 
solving $\min_{\vy = A \vc} \norm{\vc}_0$ 
is NP-hard \citep{juditsky2020statistical}, and we typically 
deal with data that is contaminated with noise, we tackle instead a relaxation of the form:
\begin{align}
	\min_{\vc \in \mathcal{C}} \norm{\vy - A \vc}^2, \label{ConvexRegressionProblem}
\end{align}
where $\mathcal{C}$ is a sparsity-inducing convex set; a typical choice is the $\ell_1$ norm ball. One of the
benefits of solving Problem~\eqref{ConvexRegressionProblem}, apart from the fact that
it can be approached with the tools of convex optimization, is that there exists a rich literature on the guarantees 
that we can expect from this recovery (see e.g., \cite{donoho2006compressed, candes2008introduction, candes2006stable, tropp2006just}).

We will present a finite-dimensional example where the basis functions are monomials 
of the features of the vector $\vx \in \mathbb{R}^{15}$ of maximum degree $d =4$, that is $f_i(\vx) = \prod_{j=1}^n x_{j}^{a_{j}}$ with $a_{j} \in \mathbb{N}$ and $\sum_{j=1}^n a_{j} \leq d$.
We generate a random vector $\vc$ that will have $5 \%$ of its entries drawn uniformly between $0$ and $10$, while the remaining entries are zero. In order to evaluate the polynomial, we generate a total of $1000$ data points $\{\vx_i \}_{i=1}^N$ from the standard multivariate Gaussian in $\mathbb{R}^{15}$, with which we will compute the output variables $\{y_i \}_{i=1}^N$. Before evaluating the polynomial these points will be contaminated with noise drawn from a standard multivariate Gaussian. For such a low number of features in the extended space, the polynomial features
can be precomputed, thus reducing the computational burden of each function
and gradient evaluation.
We leverage \texttt{MultivariatePolynomials.jl} \citep{benoit_legat_2021_4656033}
to create the input polynomial and evaluate it on the training and test data.

The sparsity-inducing convex feasible set we minimize over is the $\ell_1$ norm ball.
Solving a linear minimization problem over this feasible region generates points with only one non-zero element.
Moreover, there is a closed-form solution for these minimizers. We run the Lazy Away-Step Frank Wolfe (L-AFW) and BCG algorithms with the adaptive line search strategy from \cite{pedregosa2018step}, and compare them to Projected Gradient Descent using a smoothness estimate. We will evaluate the output solution on test points drawn in a similar manner as the training points.
The radius of the $\ell_1$ norm ball that we will use to regularize the problem will be equal to $0.95 \norm{\vc}_1$.

\begin{figure}[th!]
\centering
\includegraphics[width=0.8\linewidth]{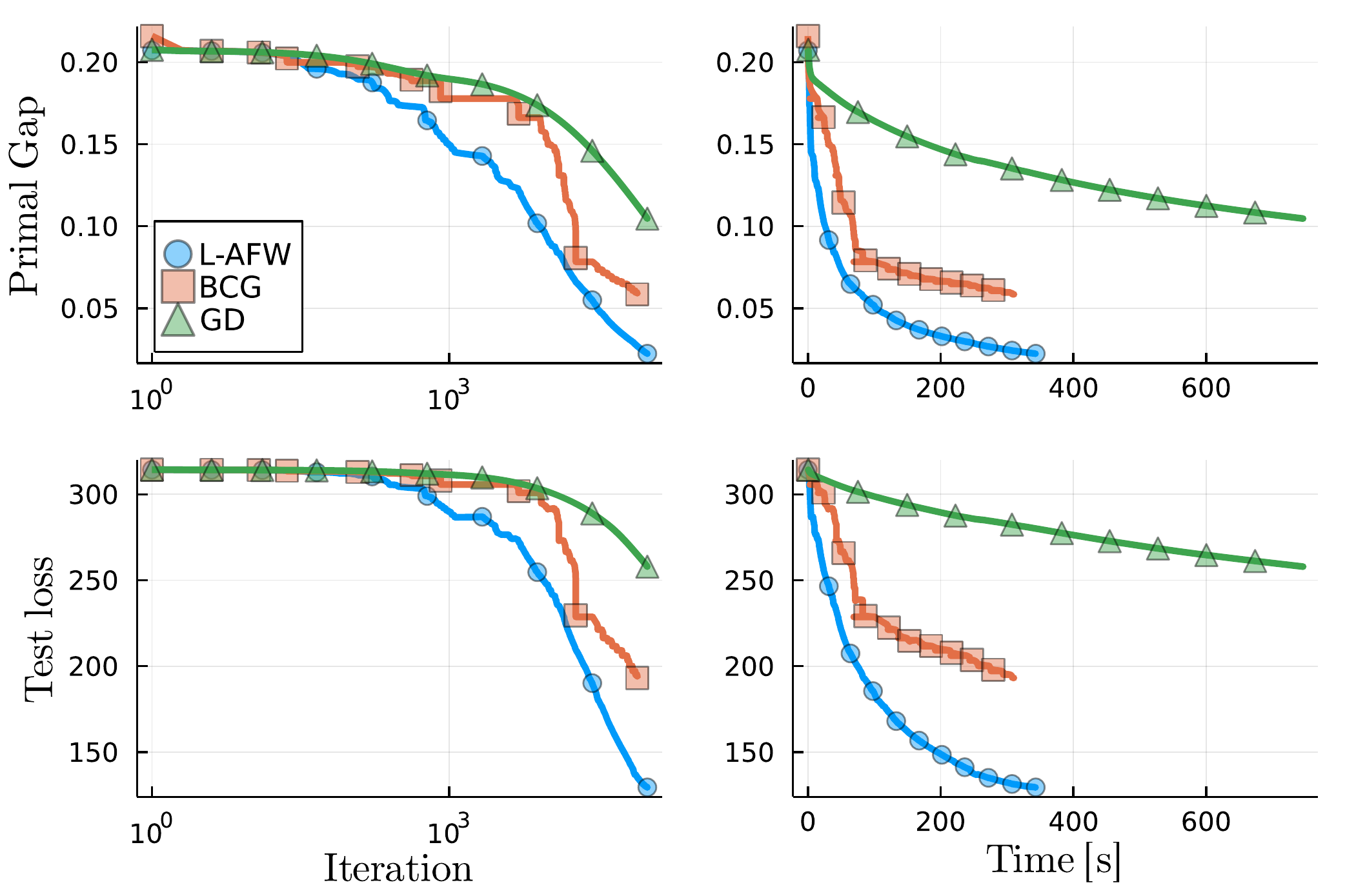}
\caption{Polynomial regression results.}
\label{fig:polyresult}
\end{figure}

The primal gap and the test error evolution during the run of the algorithms
are presented in \cref{fig:polyresult}, in terms of iteration count and wall-clock time. 
The L-AFW and BCG algorithms decrease the primal gap faster than the GD 
algorithm in terms of iteration count, due to the effectiveness of the 
adaptive stepsize strategy. This faster convergence in iteration count 
is amplified with respect to time, making the L-AFW and BCG algorithms 
converge much faster than GD, due to the fact that the L-AFW and BCG 
algorithms exploit the inherent sparsity of the extreme points of the $\ell_1$ 
ball, making the computational cost of each iteration much lower than the 
cost of each GD step. Note that none of the algorithms manage to perfectly
 recover the support of the exact coefficients of $\vc$, due to the 
 magnitude of the noise introduced in the example. This is expected, 
 given the small values that some of these coefficients have.

\subsection{Matrix completion}

Given a set of observed entries from a matrix
$Y \in \mathbb{R}^{m\times n}$, we want to compute 
a matrix $X \in \mathbb{R}^{m\times n}$ that minimizes
the sum of squared errors on the observed entries. As it stands 
this problem formulation is not well-defined or useful, as one could minimize the objective function simply by setting the observed entries of $X$ to match
those of $Y$, and setting the remaining entries of 
$X$ arbitrarily. A common way 
to solve this problem is to reduce the degrees of freedom of the problem by assuming 
that the matrix $Y$ has low rank \citep{candes2009exact, candes2010power, candes2010matrix, udell2019big}.
Finding the matrix $X \in \mathbb{R}^{m\times n}$ with minimum rank whose observed entries are equal to those of $Y$ is a non-convex problem that is $\exists \mathbb{R}$-hard \citep{bertsimas2020mixed}.
A common convex proxy for rank constraints is the use of constraints on the nuclear norm
of a matrix \citep{fazel2002matrix}, and so attempt to solve:
\begin{align}
\min_{\norm{X}_{*} \leq \tau} \sum_{(i,j)\in \mathcal{I}} \left( X_{i,j} - Y_{i,j}\right)^2, \label{Prob:matrix_completion}
\end{align}
where $\tau>0$ and $\mathcal{I}$ denotes the indices of the observed entries of $Y$. In this example, we compare the Frank-Wolfe implementation from the package
with a Projected Gradient Descent (PGD) algorithm which, after each gradient descent step, projects the iterates back onto the nuclear norm ball.
We use one of the Movielens datasets \citep{harper2015movielens} to compare the two methods, and for convenience we set $\tau$ to be $10$ times the largest singular value of the matrix $Y$ with all the available entries.
The code required to reproduce the full example is presented in the documentation.

\begin{figure}[H]
\centering
  \includegraphics[width=0.8\linewidth]{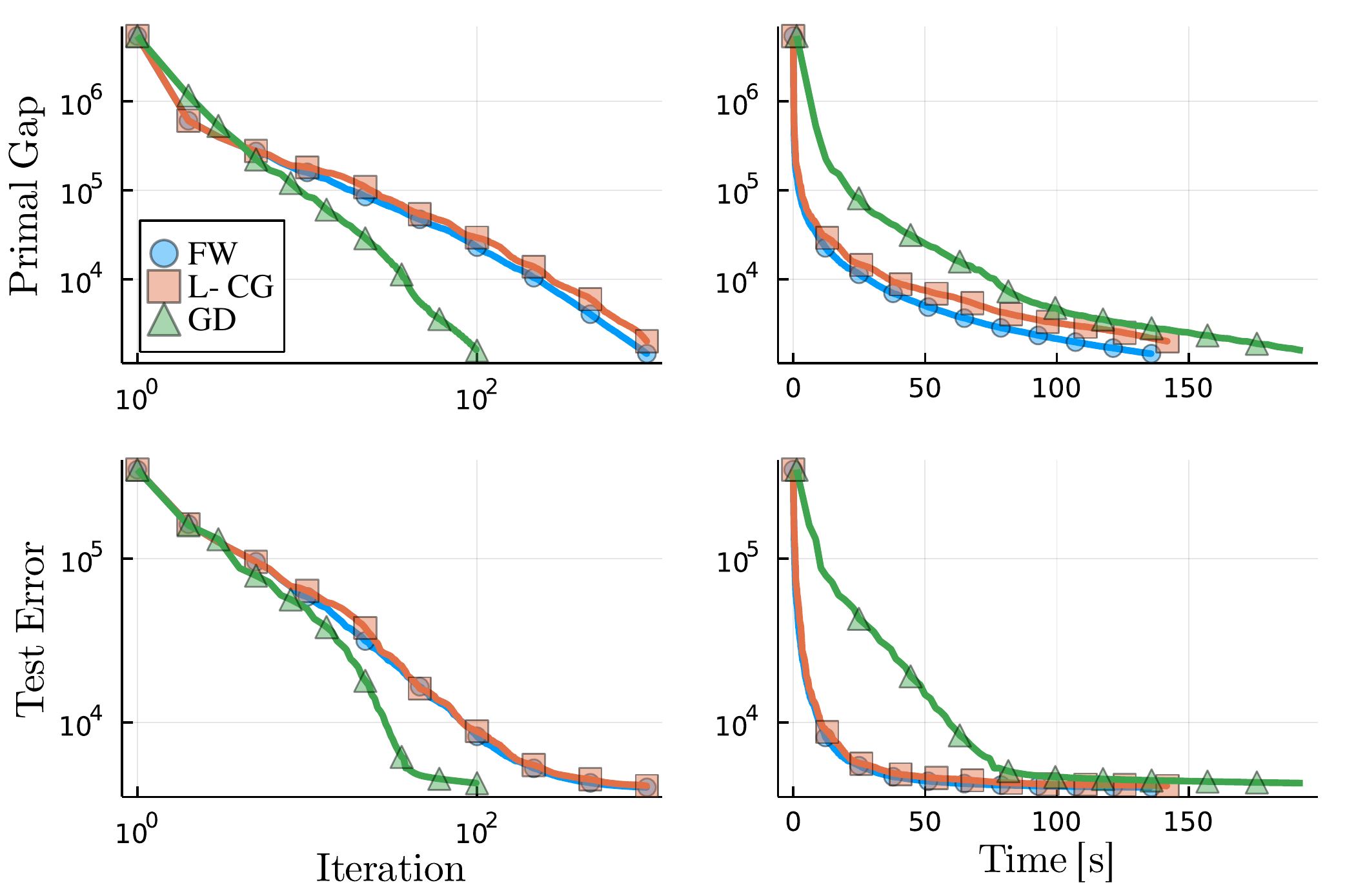}
  \caption{Comparison of standard and lazified FW with projected gradient descent on the Movielens example.}
  \label{fig:movielens}
\end{figure}

The results are presented in \cref{fig:movielens}. We can clearly observe that the computational cost
of a single PGD iteration is much higher than the cost of a FW variant step. The FW variants tested 
complete $10^3$ iterations in around $120$ seconds, while the PGD algorithm only completes $10^2$ iterations 
in a similar time frame. We also observe that the progress per iteration made by each 
projection-free variant is smaller than the progress made by PGD, as expected.
Note that, minimizing a linear function over the nuclear norm ball, in order to compute the LMO, amounts to computing the left and right singular vectors associated with the largest singular value, which we do using the \texttt{ARPACK} \citep{lehoucq1998arpack} Julia wrapper in the current example.
On the other hand, projecting onto the nuclear norm ball requires computing a full singular value decomposition. The underlying linear solver can be switched by users developing their own LMO.\\

The top two figures in \cref{fig:movielens} present the primal gap of Problem~\eqref{Prob:matrix_completion} in terms of iteration count and wall-clock time.
The two bottom figures show the performance on a test set of entries. Note that the test error stagnates for all methods, as expected.
Even though the training error decreases linearly for PGD for all iterations,
the test error stagnates quickly. The final test error of PGD is about $6\%$ higher
than the final test error of the standard FW algorithm, which is also $2\%$ smaller than the final test error of the lazy FW algorithm. We would like to stress though that the intention here is primarily to showcase the algorithms and the results are considered to be illustrative in nature only rather than a proper evaluation with correct hyper-parameter tuning. 

Another key aspect of FW algorithms is the sparsity of the provided solutions.
Sparsity in this context refers to a matrix being low-rank.
Although each solution is a dense matrix in terms of non-zeros,
it can be decomposed as a sum of a small number of rank-one terms,
each represented as a pair of left and right vectors.
At each iteration, FW algorithms add at most one rank-one term to the iterate,
thus resulting in a low-rank solution by design.
In our example here, the final FW solution is of rank at most $95$ while the lazified version provides a sparser solution of rank at most $80$. The lower rank of the lazified FW is due to the fact that this algorithm sometimes avoids calling 
the LMO if there already exists an atom (here rank-1 factor) in the cache that guarantees enough progress; the higher sparsity might help with interpretability and robustness to noise. In contrast, the solution computed by PGD is of full column rank and even after truncating the spectrum, removing factors with small singular values, it is still of much higher rank than the FW solutions. 

\subsection{Exact optimization with rational arithmetic}

The package allows for exact optimization with rational arithmetic.
For this, it suffices to set up the LMO to be rational and choose
an appropriate step-size rule as detailed below.
For the LMOs included in the package, this simply means initializing the
radius with an element type that can be promoted to a
rational\footnote{We refer the interested reader to the Julia documentation on promotion mechanisms \url{https://docs.julialang.org/en/v1/manual/conversion-and-promotion}},
for instance the integer \jlinl{1},
rather than a floating-point number such as \jlinl{1.0}.
Given that numerators and denominators can become quite large in rational arithmetic,
it is strongly advised to base the used rationals on extended-precision integer types such as \jlinl{BigInt}, i.e., we use \jlinl{Rational\{BigInt\}}.
For the probability simplex LMO with a rational radius of \jlinl{1}, the LMO would be created as follows:

\begin{jllisting}
lmo = FrankWolfe.ProbabilitySimplexOracle{Rational{BigInt}}(1)
\end{jllisting}

As mentioned before, the second requirement ensuring that the computation runs in rational
arithmetic is a rational-compatible step-size rule.
The most basic step-size rule compatible with rational optimization is the \jlinl{agnostic} step-size\
rule with $\gamma_t = 2/(2+t)$.
With this step-size rule, the gradient does not even need to be rational as long as the atom computed by the LMO
is of a rational type.
Assuming these requirements are met, all iterates and the computed solution will then be rational:

\begin{jllisting}
  n = 100
  x = fill(big(1)//100, n)
  # equivalent to { 1/100 }^100
\end{jllisting}
  
Another possible step-size rule is \jlinl{rationalshortstep} which computes the step size
by minimizing the smoothness inequality as
$\gamma_t = \frac{\langle \nabla f(\vx_t), \vx_t - \mathbf{v}_t\rangle}{2 L \|\vx_t - \mathbf{v}_t\|^2}$.
However, as this step size depends on an upper bound on the Lipschitz constant $L$ as well as the
inner product with the gradient $\nabla f(\vx_t)$, both have to be of a rational type.

\subsection{Formulating the LMO with MathOptInterface}

In this example, we project a random point onto the scaled probability simplex with the basic Frank-Wolfe algorithm
using either the specialized LMO defined in the package or a generic LP formulation using
\texttt{MathOptInterface.jl} (MOI) and \texttt{GLPK} as the underlying LP solver.
\jlinputlisting{examples/moi_optimizer.jl}

\begin{figure}[H]
\centering
  \includegraphics[width=0.8\linewidth]{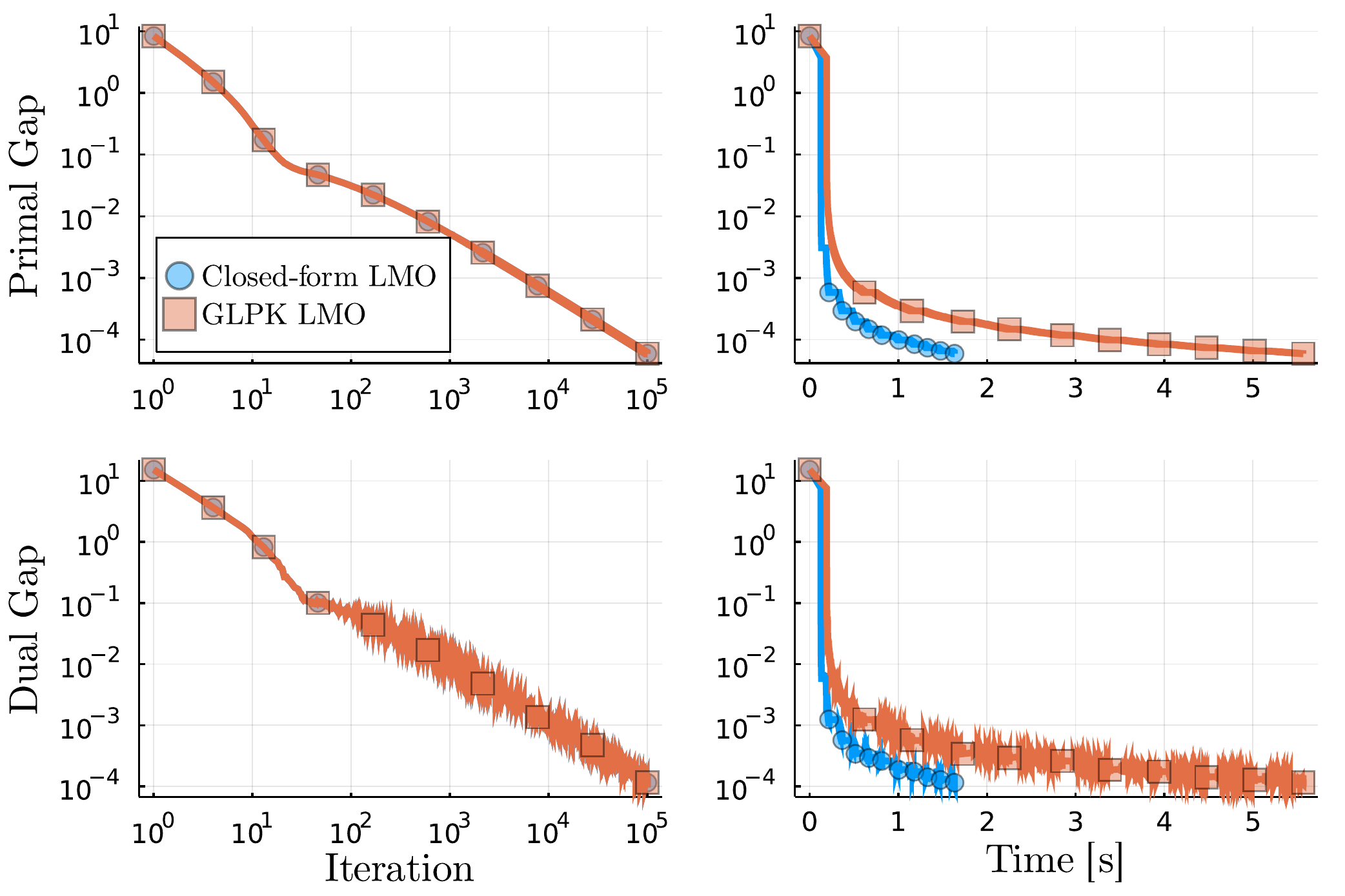}
  \caption{Performance comparison of custom LMO and MOI LP definition.}
  \label{fig:moi_vs_simplex}
\end{figure}

\noindent
The resulting primal and dual progress are presented in \cref{fig:moi_vs_simplex}.
Since we use identical algorithmic parameters, the per-iteration
progress is identical for the two LMOs. On a per iteration basis, the computational cost
of calling the closed form LMO is lower, as we only allocate an output vector, and we 
avoid having to solve a generic linear optimization problem.

\subsection{Doubly stochastic matrices}

\begin{figure}[H]
\centering
  \includegraphics[width=0.8\linewidth]{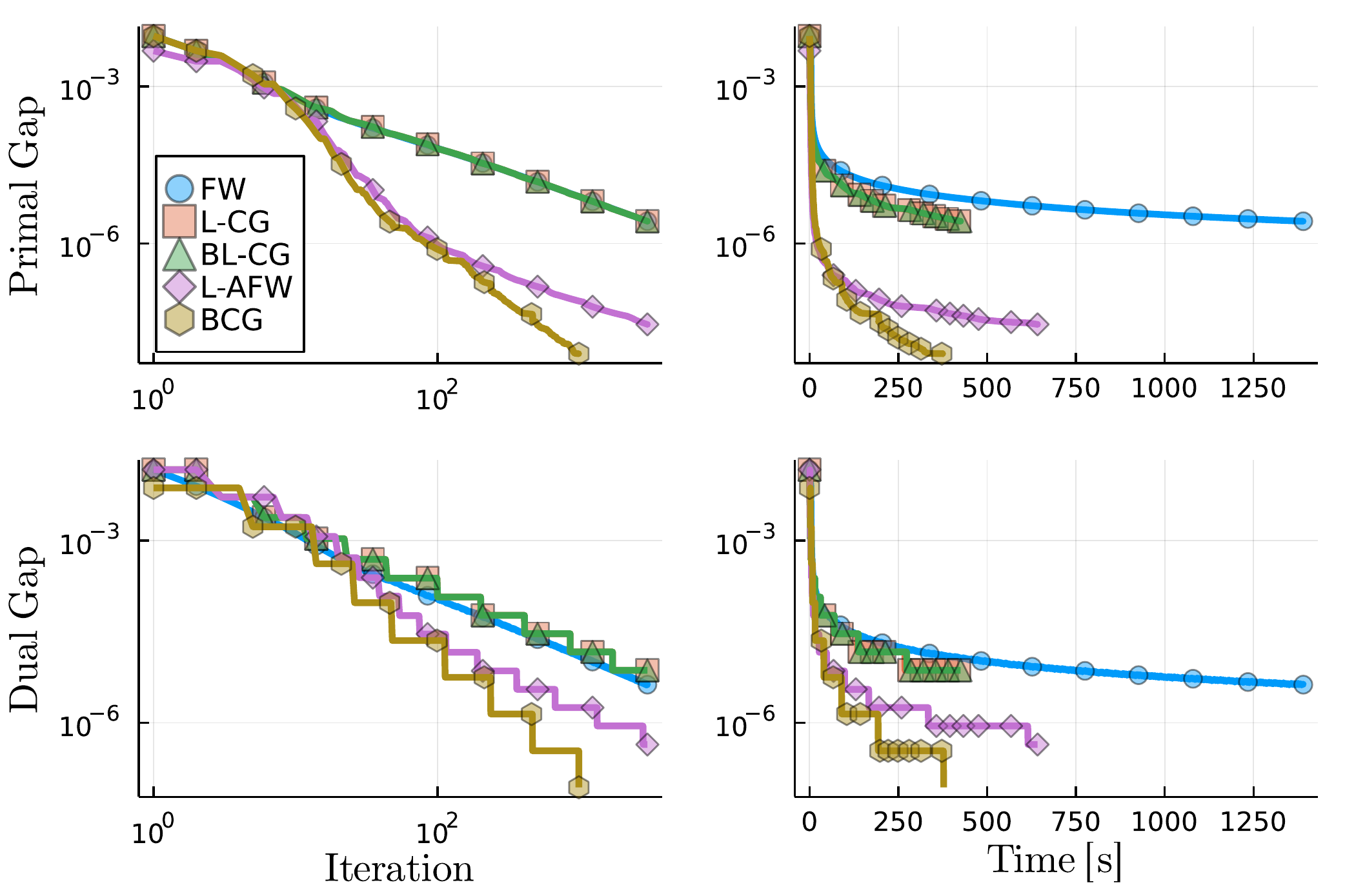}
  \caption{Lazified and eager variants on the Birkhoff polytope.
  FW stands for the classic Frank-Wolfe algorithm, L-CG uses a caching LMO with an unbounded cache size while BL-CG uses a bounded cache with $500$ elements, L-AFW is the lazified AFW algorithm, and BCG is the Blended Conditional Gradient algorithm.}
  \label{fig:lazy_birkhoff}
\end{figure}

The set of doubly stochastic matrices or Birkhoff polytope appears in various combinatorial
problems including matching and ranking.
It is the convex hull of permutation matrices, a property of interest for FW algorithms
because the individual atoms returned by the LMO only have $n$ non-zero entries for matrices in
$\mathbb{R}^{n\times n}$.
A linear function can be minimized over the Birkhoff polytope using the Hungarian algorithm.
This LMO is substantially more expensive than minimizing a linear function over the $\ell_1$-ball norm,
and thus the algorithm performance benefits from lazification.
We present the performance profile of several FW variants in the following example
for $n = 200$. The results are presented in \cref{fig:lazy_birkhoff}.

The per-iteration primal value evolution is nearly identical for FW and the lazy cache variants.
We can observe a slower decrease rate in the first $10$ iterations of BCG for both the primal value
and the dual gap. This initial overhead is however compensated after the first iterations,
BCG is the only algorithm terminating with the desired dual gap of $10^{-7}$ and not with the
iteration limit. In terms of runtime, all lazified variants outperform the standard FW,
the overhead of allocating and managing the cache are compensated by the reduced number of calls to the LMO.

\section{Final Comments and Future Work}\label{sec:futurework}

The \href{https://github.com/ZIB-IOL/FrankWolfe.jl}{\texttt{FrankWolfe.jl}} package will be further extended over time and we welcome contributions, reporting of issues and bugs, as well as pull requests under the package \href{https://github.com/ZIB-IOL/FrankWolfe.jl}{GitHub repository}.
A few prominent features to come in the near future will be:

\begin{enumerate}
	\item Stronger interfacing with the broader Julia optimization ecosystem.
	\item Interfacing with Python similar to, e.g., \href{https://github.com/SciML/diffeqpy}{diffeqpy}.
	\item Implementating closely related variants such as matching pursuit variants.
\end{enumerate}

\section*{Acknowledgements}

Research reported in this paper was partially supported through the Research Campus Modal funded by the German Federal Ministry of Education and Research (fund numbers 05M14ZAM,05M20ZBM) and the Deutsche Forschungsgemeinschaft (DFG) through the DFG Cluster of Excellence MATH+.

\bibliography{refs}
\bibliographystyle{icml2021}

\appendix \newpage

\section{Full nuclear norm example}\label{app:movielens}

Setup and data preprocessing:
\jlinputlisting{examples/movielens.jl}

\newpage
PGD and FW run:
\jlinputlisting{examples/movielens2.jl}

\end{document}